\documentclass[11pt]{article}
\usepackage{amssymb}
\newtheorem{theorem}{Theorem}

\newtheorem{remark}{Remark}

\def\defi{\stackrel{{\scriptscriptstyle \Delta}}{=}}

\def\d{\delta}
\def\o{\omega}

\def\F{{\cal F}}

\def\Ind{{\mathbb{I}}}

\def\Arg{{\rm Arg\,}}

\def\Re{{\rm Re\,}}

\def\R{{\bf R}}
\def\E{{\bf E}}

\def\H{{\cal H}}

\def\L{L}

\def\s{\delta}

\def\C{{\bf C}}
\def\W{{\cal W}^*}
\def\W{{\cal W}}

\def\oo{\bar}
\def\s{\sigma}

\def\p{\partial}
\def\G{\Gamma}

\def\M{{\cal M}}

\def\L{{\cal L}}

\newcommand{\be}{\begin{equation}}
\newcommand{\ee}{\end{equation}}
\newcommand{\bd}{\begin{displaymath}}
\newcommand{\ed}{\end{displaymath}}
\newcommand{\ba}{\begin{array}{ll}}
\newcommand{\ea}{\end{array}}
\newcommand{\baa}{\begin{eqnarray}}
\newcommand{\eaa}{\end{eqnarray}}
\newcommand{\baaa}{\begin{eqnarray*}}
\newcommand{\eaaa}{\end{eqnarray*}}
\font\sm=cmr10

\def\oo{\bar}

\title{Regularity of a inverse problem for generic parabolic equations }
\author{
Nikolai Dokuchaev\\ {\sm  Department of Mathematics, Trent
University, Ontario, Canada}}
 
\begin{document}
 \vspace{-0.5cm}      \maketitle
\begin{abstract}
The paper studies some inverse boundary value problem for simplest
parabolic equations such that  the homogenuous Cauchy condition is
ill posed at initial time. Some regularity of the solution is
established  for a wide class of boundary value inputs.
\end{abstract}
Parabolic equations such as heat equations have fundamental
significance for natural sciences, and various boundary value
problems for them were widely studied including well-posed
problems as well as the so-called inverse and ill-posed problems
that often are very significant for applications  (see, e.g., Beck
(1985)). The present paper investigates a inverse boundary value
problem on semi-plane for homogenuous parabolic equations with
homogenuous Cauchy condition at initial time and with Dirichlet
condition on the boundary of semi-plane. The parabolic equation is
the equation of a backward type that is usually solvable with
Cauchy condition at terminal time. However, we consider this
equation with Cauchy condition at initial time so that the problem
is an inverse problem. One may think that these problems are
always ill-posed in the sense that there is no regularity of
solutions such as prior estimates for the solution via some norms
of free terms  (see, e.g., Beck (1985), Tikhonov and Arsenin
(1977)). We found a class of very generic parabolic equation with
constant coefficients and with certain sign of the drift
coefficient (i.e., the coefficient for first derivative) such that
the inverse problem has some regularity. More precisely, we  found
a wide enough class of inputs in the boundary value condition on
the boundary of the semi-plane that ensures regularity in a form
of prior energy type estimates. This class of inputs is everywhere
dense in the class of $L_2$-integrable functions; it includes
differentiable functions.
\section{The problem setting}
Let us consider the following boundary value problem on
semi-plane: \baa &&a\frac{\p u}{\p t}(x,t)+\frac{\p^2 u}{\p
x^2}(x,t)+b\frac{\p u}{\p x}(x,t)+c u(x,t)=0,\quad x>0,\
t>0,\nonumber\\ &&u(x,0)\equiv 0,\quad x>0, \nonumber\\
&&k_0u(0,t)+k_1\frac{\p u}{\p x}(0,t)\equiv g(t),\quad t>0.
\label{parab} \eaa
 Here $x>0$, $t>0$, and $a,b,c, k_0,k_1\in\R$ are constants.
 \par
  We assume that \baa a>0,\quad
 b>0,\qquad k_0^2+k_1^2>0,\quad k_0k_1\le 0.
 \label{kk}\eaa
 The assumption that $a>0$ and the presence of the initial condition at $t=0$
 makes problem (\ref{parab}) an inverse problem (see, e.g.,
Beck (1985), Tikhonov and Arsenin (1977)).
\par
We assume that $u(x,t)\equiv 0$ and $g(t)\equiv 0$ for $t<0$.
\par
 Let
$\G$ denotes the set of all functions $g:\R\to\R$ such that
$g(t)=0$ for $t<0$ and with finite norm
\baaa\|g\|_{W^1_2(\R)}\defi\left\|g\right\|_{L_2(\R)}+
\left\|\frac{\p g}{\p t}\right\|_{L_2(\R)}. \eaaa
\begin{remark}  Functions $g\in\G$ are continuous on $\R$ and vanising on $t<0$ since $dg(t)/dt\in
L_2(\R)$. For instance, $g(t)=e^{-t}\sin t\in \G$, but
$g(t)=e^{-t}\cos t\notin\G$.
\end{remark}
\par
Let $D\defi \R\times\R^+$. Let  $\W$ be the space of the functions
$v=v(x,t):\R\times\R^+\to\R$ such that $v(x,t)\equiv 0$ for $t<0$
and with finite norm \baaa \|v\|_{\W}
\defi \|v\|_{L_2(D)}+\Bigl\|\frac{\p u}{\p x}
\Bigr\|_{L_2(D)}+\Bigl\|\frac{\p^2 v}{\p x^2}
\Bigr\|_{L_2(D)}+\Bigl\|\frac{\p v}{\p t}\Bigr\|_{L_2(D)}. \eaaa
The class $\W$ is such that all the equations presented in problem
(\ref{parab}) are well defined for any $u\in\W$. Let us show this.
If $v\in \W$, then, for any $t_*>0$, we have that
$v|_{\R^+\times[0,t_*]}\in C([0,t_*],L_2(\R^+))$ as a function of
$t\in[0,t_*]$. Hence the initial condition at time $t=0$ is well
defined as an equality in $L_2(\R^+)$. Further, for any $x_*>0$,
we have that  $v|_{[0,x_*]\times\R^+}\in C([0,x_*],L_2(\R^+))$ and
$\frac{\p v}{\p x}\Bigr|_{[0,x_*]\times\R^+}\in
C([0,x_*],L_2(\R^+))$ as functions of $x\in[0,x_*]$. Hence the
functions $g_0(t)\defi v(0,t)$, $g_1(t)\defi
\frac{du}{dx}(x,t)|_{x=0}$ are well defined as elements of
$L_2(\R^+)$, and the boundary value condition at $x=0$ is well
defined as an equality in $L_2(\R^+)$.
\section{The main
result}
\begin{theorem}\label{ThM} Let condition (\ref{kk}) be satisfied, and let  \baa
 \mu\defi b^2/4<c.
 \label{mu}
 \eaa
Then there exists a unique solution $u(x,t)$   in the class
 $\W$
 of problem
 (\ref{parab}) in the domain $D$ for any $g\in\G$. Moreover, there exists a constant $C=C(a,b,c,k_0,k_1)$ such
 that
\baa   \|u\|_{\W}\le C \left\|g\right\|_{W^1_2(\R)}. \label{estM}
\eaa
\end{theorem}
\begin{remark}
It will be seen from the proof that it is crucial that
$u(x,0)\equiv 0$ in the Cauchy condition and that the parabolic
equation is homogenuous. We cannot  extend the result for
 non-zero initial conditions or non-zero
free term in the parabolic equation.
\end{remark}
\par
 The following
theorem shows that  assumption (\ref{mu}) is not really
restrictive
 if we are not interested  in the properties of the solutions for
 $T\to +\infty$, as can happen if we deal with solutions  a finite time interval.
\begin{theorem}\label{ThM2}
Let condition (\ref{kk}) holds, but (\ref{mu}) does not hold. Let
$M$ be such that $g(t)e^{-Mt}\in \G$ and $b^2/4<c+M$. Then problem
(\ref{parab}) has a unique solution $u$ such that $u_M\in \M$,
where $u_M(x,t)\defi e^{-Mt}u(x,t)$.
\end{theorem}
\par
{\it Proof of Theorem \ref{ThM}}. Let $\R^+\defi[0,+\infty)$,
$\C^+\defi\{z\in\C:\ \Re z>  0\}$. For $v\in L_2(\R)$, we denote
by $\F v$ and $\L v$ the Fourier and the Laplace  transforms
respectively
 \baa\label{U}
V(i\o)=(\F v)(i\o)\defi\frac{1}{\sqrt{2\pi}}\int_{\R}e^{-i\o
t}v(t)dt, \quad\o\in\R,\eaa
 \baa\label{Up}
V(p)=(\L v)(p)\defi\frac{1}{\sqrt{2\pi}}\int_{0}^{\infty}e^{-p
t}v(t)dt, \quad p\in\C^+. \eaa
\par
Let  $H^r$ be the Hardy space of holomorphic on $\C^+$ functions
$h(p)$ with finite norm
$\|h\|_{H^r}=\sup_{k>0}\|h(k+i\o)\|_{L^r(\R)}$, $r\in[1,=\infty]$
(see, e.g., Duren (1970)).
\par
Let $u\in\W$ be a solution of (\ref{parab}). Set  $g_0(t)\defi
u(0,t)$, $g_1(t)\defi \frac{du}{dx}(x,t)|_{x=0}$. As was discussed
above, the functions $g_k$ are well defined as elements of
$L_2(\R^+)$.
\par
Let $G\defi\L g$. We have that $G\in H^2$. Let $G_k(p)\defi \L
g_k$ and $U\defi\L u$ be defined for $p\in\C^+$. They are well
defined since $u\in\W$; in addition,  $G_k\in H^2$.
\par
For functions $V:\R^+\times\oo\C^+\to\C$, where $\oo\C^+=\{z:\ \Re
z\ge 0\}$, we introduce semi-norms \baaa &&\|V\|_{L_{22}}
\defi \biggl(\int_{\R^+}dx \int_\R|V(x,i\o|^2d\o\biggr)^{1/2},\\
&& \|V\|_{\W^*}
\defi \|V\|_{L_{22}}+\Bigl\|\frac{\p V}{\p
x}\Bigr\|_{L_{22}}+\Bigl\|\frac{\p^2 V}{\p x^2}
\Bigr\|_{L_{22}}+\|\frac{\p V}{\p t} \Bigr\|_{L_{22}}. \eaaa and
norms  \baaa &&\|V\|_{L^H_{22}}
\defi \biggl(\int_{\R^+} \|V(x,\cdot)\|^2_{H^2}dx\biggr)^{1/2},\\
&&\|V\|_{\H} \defi \|V\|_{L^H_{22}}+\Bigl\|\frac{\p V}{\p
x}\Bigr\|_{L^H_{22}}+\Bigl\|\frac{\p^2 V}{\p x^2}
\Bigr\|_{L^H_{22}}+\Bigl\|\frac{\p V}{\p t} \Bigr\|_{L^H_{22}}.
\eaaa
 Instead of (\ref{parab}), consider the problem \baa &&apU(x,p)+\frac{\p^2
U}{\p x^2}(x,p)+b\frac{\p U}{\p x}(x,p)+c U(x,t)=0,\quad x>0,
\nonumber\\ &&k_0U(0,p)+\frac{\p U}{\p x}(0,p)\equiv G(p),\quad
p\in \C^+ \label{parabU}\eaa subject to the following condition
 \baa U(x,\cdot),\frac{\p
U}{\p x}(x,\cdot), \frac{\p^2 U}{\p x^2}(x,\cdot)\in
H^2\quad\hbox{for a.e}\quad x>0 ,\qquad \|U\|_{\H}<+\infty.
 \label{L2}
 \eaa
 \par
 Let
$\lambda_k=\lambda_k(p)$ be the roots of the equation
$\lambda^2+b\lambda +(c+ap)=0$ defined for $p\in\C^+$ as
 $\lambda_{1}\defi -b/2-
\sqrt{\mu-ap}$ and  $\lambda_{2}\defi -b/2+ \sqrt{\mu-ap}$, where
$\mu= b^2/4-c<0$. We mean the branch of the square root such that
$\Arg\sqrt{\mu-ap}\in
[-\pi/2,+\pi/2]$ 
and $\Re\sqrt{\mu-ap}\ge 0$. Under these assumptions, the function
$\sqrt{\mu-ap}$ is holomorphic and does not have zeros in $\C^+$.
We have that \baa &&\Re \lambda_1(p)\le -\frac{b}{2}, \quad
p\in\C_+,\nonumber\\ &&\exists \d>0,\ \o_*>0:\
\Re\lambda_2(i\o)>\d \quad\hbox{if}\quad \o\in\R,\
|\o|\ge\o_*.\hphantom{xxx}
 \label{relam}
\eaa
\par In addition, we have that the functions $\lambda_k(p)$ are
holomorphic in $\C^+$, and \baa
(\lambda_1(p)-\lambda_2(p))^{-1}\in H^{\infty},\quad
\lambda_k(p)(\lambda_1(p)-\lambda_2(p))^{-1}\in H^{\infty},\quad
k=1,2, \nonumber\\ (k_0+k_1\lambda_1(p))^{-1}\in H^{\infty},\quad
\lambda_1(k_0+k_1\lambda_1(p))^{-1}\in H^{\infty}. \label{lamH}
 \eaa
 The last two statements here follow from (\ref{kk}).
 Let $$N\defi
\left\|\frac{1}{\lambda_1-\lambda_2}\right\|_{H^\infty}+\sum_{k=1,2}
\left\|\frac{\lambda_k}{\lambda_1-\lambda_2}\right\|_{H^\infty}
+\left\|\frac{1}{k_0+k_1\lambda_1}\right\|_{H^\infty}
+\left\|\frac{\lambda_1}{k_0+k_1\lambda_1}\right\|_{H^\infty}.
 $$
It can be seen also that the functions $ e^{x\lambda_k(p)}$ are
holomorphic in $\C^+$ for any $x>0$.
 \par
  For any $x
>0$, the unique solution of (\ref{parabU}) is
 \baaa
 U(x,p)&=&
\frac{1}{\lambda_1-\lambda_2}\biggl((G_1(p)-\lambda_2G_0(p))e^{\lambda_1
x}+(G_1(p)-\lambda_1G_0(p))e^{\lambda_2x}\biggl).
 \eaaa
  This can be derived, for instance,  using
 Laplace transform method applied to linear ordinary differential equation
 (\ref{parabU}), and having in mind that
   \baaa \frac{1}{\lambda^2+b\lambda +c-ap}= \frac{1}{(\lambda-\lambda_1)(\lambda-\lambda_2)}=
\frac{1}{\lambda_1-\lambda_2}\left(\frac{1}{\lambda-\lambda_1}
-\frac{1}{\lambda-\lambda_2}\right),\\
\frac{\lambda}{\lambda^2+b\lambda
+c-ap}=\frac{\lambda}{(\lambda-\lambda_1)(\lambda-\lambda_2)}=
\frac{1}{\lambda_1-\lambda_2}\left(\frac{\lambda_1}{\lambda-\lambda_1}
-\frac{\lambda_2}{\lambda-\lambda_2}\right). \eaaa
\par
 Let
represent $U$  as
 $U(x,p)=U_1(x,p)+U_2(x,p)$, where
 \baaa
&&U_1(x,p)=e^{\lambda_1 x}J_1(p),\qquad J_1(p)=
\frac{1}{\lambda_1-\lambda_2}(G_1(p)-\lambda_2G_0(p)),\\
&&U_2(x,p)=e^{\lambda_2x}J_2(p),\qquad J_2(p)
=\frac{1}{\lambda_1-\lambda_2}(G_1(p)-\lambda_1G_0(p)).
 \eaaa
 By (\ref{relam}), the fact that $u\in \W$ implies that (\ref{L2}) holds. Let us show that
  \baa
 \|U_1\|_{L_{22}}<+\infty.\label{L22}
 \eaa
By (\ref{relam}), $|e^{x\lambda_1(p))}|\le e^{-bx/2}<1$,
$p\in\C^+$.
  It
follows that \baaa  \|e^{x\lambda_1(p)} J_1(p)\|_{H^2}\le
\sup_{p\in\C^+}|e^{x\lambda_1(p))}|\|J_1\|_{H^2}\le e^{-bx/2 }\|
J_1\|_{H^2}\le Ne^{-bx/2 }\sum_{k=0,1}\| G_k\|_{H^2}. \eaaa
 Then (\ref{L2}) and (\ref{L22}) imply that
 \baa\|U_2\|_{L_{22}}<+\infty.
 \label{U2}\eaa
 Further,
 \baaa
+\infty > \|U_2\|^2_{L_{22}}=
\int_{\R^+}dx\int_{\R}|U_2(x,i\o)|^2d\o=
\int_{\R^+}dx\int_{\R}|e^{\lambda_2(i\o)x}J_2(i\o)|^2
d\o\\=\int_{\R^+}dx  \int_{\R}e^{\Re
\lambda_2(i\o)x}|J_2(i\o)|^2d\o\ge \int_{\R^+}dx e^{\d x}
\int_{\o: |\o|\ge \o_*}|J_2(i\o)|^2d\o.
 \eaaa
 Note that the  $J_2(i\o)$ is vanishing on $\{\o: |\o|\ge
 \o_*\}$. Since $J_2\in H^2$, it follows that
 $$
 J_2=(G_1(p)-\lambda_1G_0(p))e^{\lambda_2x}\equiv 0,$$
 i.e.,
 $$
 G_1(p)=\lambda_1G_0(p).$$
 Remind that $k_0G_0(p)+k_1G_1(p)=G(p)$. It follows that
\baaa &&k_0G_0(p)+k_1\lambda_1G_0(p)=G(p),\qquad
 G_0(p)=(k_0+k_1\lambda_1(p))^{-1}G(p),
\\
&&
J_1(p)= G_0(p),\qquad U(x,p)=U_1(x,p)=e^{\lambda_1 x}G_0(p).
\eaaa
\par
Let us estimate $\|U\|_{\H}$.
\par
By (\ref{relam}), $|e^{x\lambda_1(p))}|\le e^{-bx/2}<1$.
  It
follows that \baaa  &&\|p^me^{x\lambda_1(p)} G_0(p)\|_{H^2}\le
e^{-bx/2}\|p^mG_0(p)\|_{H^2}\\ &&\le
e^{-bx/2}\|(k_0+k_1\lambda_1(p))^{-1}\|_{H^\infty}\|p^mG\|_{H^2}\le
Ne^{-bx/2 }\|p^mG\|_{H^2},\quad m=0,1. \eaaa
\par
It  follows from the above estimate that \baa
\|p^mU\|_{L^H_{22}}\le N C_1 \left\|p^m G\right\|_{H^2},\quad
m=0,1. \label{est1} \eaa
\par
Further, we have that \baa \frac{\p U}{\p x}(x,p)=
G_0(p)\lambda_1e^{\lambda_1 x}=
\frac{\lambda_1}{k_0+k_1\lambda_1}G(p)\lambda_1e^{\lambda_1 x}.
 \label{est2}\eaa
  We obtain again
that \baa \left\|\frac{\p U}{\p x}\right\|^2_{L^H_{22}}=
\int_{\R^+}\left\|\frac{\p U}{\p x}(x,p)\right\|^2_{H^2}dx\le N
C_2
\int_{\R^+}\left\|e^{\lambda_1(p)x}G(p)\right\|^2_{H^2}dx\nonumber\\
\le N C_2 \int_{\R^+}e^{-bx/2}\left\|G(p)\right\|^2_{H^2}dx\le
C_3\|G\|_{H^2}.\label{est3}\eaa By (\ref{parabU}), $\p^2U/\p x^2$
can be expressed as a linear combination of $U$, $pU$, and $\p
U/\p x$. By (\ref{est1})-(\ref{est3}), \baaa
\int_{\R^+}\left\|\frac{\p^2 U}{\p x^2}(x,p)\right\|^2_{H^2}dx\le
C_4 \left( \int_{\R^+}\left\|\frac{\p U}{\p
x}(x,p)\right\|_{H^2}^2dx+\sum_{m=0,1}
\int_{\R^+}\left\|p^mU(x,p)\right\|_{H^2}^2dx\right). \eaaa
 It follows that \baa \int_{\R^+}\left\|\frac{\p^2 U}{\p
x^2}(x,p)\right\|_{H^2}^2dx\le C_5 (\left\|
G\right\|_{H^2}^2+\left\|pG(p)\right\|_{H^2}^2).\label{est4} \eaa
Here $C_k$ are constants that depend on $a,b,c,k_0,k_1$. By
(\ref{est1})-(\ref{est4}), estimate (\ref{L2}) holds.
\par
  Let $u(x,\cdot)\defi \F^{-1}U(x,i\o)|_{\o\in \R}$.
By (\ref{est1}), it follows that the corresponding inverse Fourier
transforms $u(x,\cdot)=\F^{-1}U(x,i\o)|_{\o\in \R}$, $\frac{\p
u}{\p t} (x,\cdot)=\F^{-1}(pU(x,i\o)|_{\o\in \R})$ are well
defined and are vanishing for $t<0$.
 In addition, we have that
$\overline{U(x,i\o)}=U(x,-i\o)$ (for instance,
 $\overline{G_k(i\o)}=G_k(-i\o)$,
 $\overline{e^{x\lambda_k(i\o)}}=e^{x\lambda_k(-i\o)}$, etc). It follows that
 the inverse of Fourier transform
 $
u(x,\cdot)=\F^{-1}U(x,\cdot)
 $
 is real. By
(\ref{L2}), estimate (\ref{estM}) holds. Therefore, $u$
   is the solution of (\ref{parab}) in $\W$.
  The uniqueness is ensured by the linearity of the problem, by
  estimate (\ref{estM}), and by the fact that $\L u(x,\cdot)$, $\L (\p^k u(x,\cdot)/\p x^k)$,
  and  $\L (\p u(x,\cdot/\p t)$ are well defined on $\C^+$
   for any function $u$ from
  $\W$, i.e, (\ref{parabU}) must be satisfied together with (\ref{L2}).
  This completes the proof of Theorem \ref{ThM}. $\Box$
\par
{\it Proof of Theorem \ref{ThM2}}.
  Rewrite
 the parabolic equation as the one with $c$ replaced by $c+M$ and $g(t)$ replaced by
 $g(t)e^{-Mt}$. By Theorem \ref{ThM},
  solution  $u_M\in\M$ of the new equation exists. Clearly,
 $u(x,t)=e^{Mt}u_M(x,t)$ is the solution of of the original problem.  $\Box$
\section{Some application}
\subsection*{Total absorbing on the boundary}
Let $T>0$ be given. Let us consider the following well-posed
boundary value problem on semi-plane: \baa &&\frac{\p v}{\p
t}(x,t)+a\frac{\p^2 v}{\p x^2}(x,t)+b\frac{\p v}{\p
x}(x,t)+cv(x,t)=0,\quad \quad x>0,\ t\in [0,T],\nonumber\\
&&v(x,T)\equiv v_*(x), \nonumber\\ &&k_0v(0,t)+k_1\frac{\p v}{\p
x}(0,t)\equiv g(t). \label{parabN} \eaa
 Here  $a,b,c, k_0,k_1\in\R$ are constants such that (\ref{kk}) holds.
  \begin{theorem}\label{cor}
  For any $g\in\G$, there exists $v_*\in L_2(\R^+)$ such that  $v(x,0)\equiv 0$,
  where  $v\in\W$ is the solution of   well-posed problem
(\ref{parabN}).
\end{theorem}
{\it Proof}. It suffices to take the solution $u\in\W$ of problem
(\ref{parab}) and take $v_*(x)\defi u(x,T)$, $v=u$. $\Box$
\par
Note that  one can rewrite problem (\ref{parabN})  as a well posed
problem for forward parabolic equation with initial time at time
$t=0$ via time change $t\to T-t$. In that case, the phenomena
described in Theorem \ref{cor} looks more impressive.
\subsection*{Restoring past distributions of diffusion processes}
Consider the following stochastic process
\be
\label{yxs} y^{x}(t)=x+bt+\s w(t). \ee Here $x\ge 0$,  $w(t)$ is a
scalar Wiener process, $b>0$ and $\s>0$ are constants.  \par
 Let $a$ be a random number such that $a\ge 0$ and it has the
probability density function $\rho\in L_2(\R^+)$ which is supposed
to be unknown. We assume also that $a$ is independent from
$w(t)-w(t_1)$ for all $t>t_1\ge 0$.
\par
 Let $y(t)=y^{a}(t)$ be the solution of  Ito equation
(\ref{yxs}) with the initial condition $y(0)=a$.
 Set
$\tau^{a}\defi\min\,\{t>0: y^{a}(t)=0\}$.
\par
Let $\Ind$ denotes the indicator function of an event.\par
 Let  $p(x,t)$ be
the  probability density function of the process $y^{a}(t)$ if
this process is killed at $0$ and inside $(0,+\infty)$ with the
rate of killing $c$ (case of $c>0$ is not excluded for the sake of
generality). More precisely, $p$ is such that $$ \int_{B} p(x,t)dx
=\E e^{ct}\Ind_{(y^a(t)\in B})\Ind_{\{ \tau^{a}\ge t\}} $$ for any
domain $B\subset \R^+$. It is known that evolution of $p$ is
described by the parabolic equation being ajoint to (\ref{parab})
with the boundary value conditions $p(x,0)\equiv\rho(x)$,
$p(0,t)\equiv 0$.
\par
\par For $g\in\G$, let  $u_g=u_g(x,t)$ be the
solution of  the problem (\ref{parab}), where $k_0=1$, $k_1=0$.
\par
Let $T>0$ be given, and let $\Psi_g(x)\defi u_g(x,T)$.
\par
\begin{theorem}\label{Thprobd} For all
functions $g\in\G$ and all $c\in\R$, \be\E
e^{c\tau^a}g(\tau^{a})\Ind_{\{ \tau^{a}<T\}}=-\E
e^{cT}\Psi_g(y(T))\Ind_{\{ \tau^{a}\ge T\}} =-\int_{\R^+}
p(x,T)\Psi_g(x)dx. \label{probd}\ee
\end{theorem}
Theorem \ref{Thprobd}  allows to solve effectively the following
inverse problem: find the distribution of $\tau^{a}\Ind_{\{
\tau^{a}<T\}}$   for unknown distribution of $a$ using the
"future" values of $p(x,T)$ only. More precisely, one can find the
values of the expectations at the left hand side of (\ref{probd})
for all $g\in\G$ using only  $p(x,T)$ via the following algorithm:
\begin{itemize}
\item[(a)] Find $u_g$ as the solution of (\ref{parab});
\item[(b)] Find $\Psi_g=u_g(\cdot,T)$;
\item[(c)] Using known $p(x,T)$, calculate the integral at the right hand side of
(\ref{probd}).
\end{itemize}
This approach does not require regularization of a ill-posed
problem such as in Beck (1985) or Tikhonov and Arsenin (1977).
\par
{\it Proof of Theorem \ref{Thprobd}}. Clearly, \baaa\E
e^{c(\tau^a\land T)} u_g(y(\tau^a\land T),\tau^a\land T)=\E
e^{c\tau^a}g(\tau^{a})\Ind_{\{ \tau^{a}<T\}}+\E
e^{cT}\Psi_g(y(T))\Ind_{\{ \tau^{a}\ge T\}}. \eaaa By Ito formula,
$$\E e^{c(\tau^a\land T)} u_g(y(\tau^a\land T),\tau^a\land T) =\E
u_g(a,0)=0. $$  In addition, \baaa \E e^{cT}\Psi_g(y(T))\Ind_{\{
\tau^{a}\ge T\}}=\int_{\R^+} p(x,T)\Psi_g(x)dx. \eaaa Then the
result follows. $\Box$
\subsection*{References}
$\phantom{xxi}$Beck, J.V. {\it Inverse Heat Conduction.} 1985.
John Wiley and Sons, Inc..
\par
Duren, P. {\it Theory of $H^p$-Spaces.} 1970. Academic Press, New
York.
\par
Tikhonov, A. N. and Arsenin, V. Y. {\it Solutions of Ill-posed
Problems.} 1977. W. H. Winston, Washington, D. C.
\end{document}